\documentstyle{article}
\def\hybrid{\topmargin 0pt      \oddsidemargin 0pt
        \headheight 0pt \headsep 0pt
        \textwidth 6.5in        % US paper
        \textheight 9.0in         % US paper
        \marginparwidth 0.0in
        \parskip 0pt plus 1pt   \jot = 1.5ex}
\catcode`\@=11
\def\marginnote#1{}

\newcount\hour
\newcount\minute
\newtoks\amorpm
\hour=\time\divide\hour by60
\minute=\time{\multiply\hour by60 \global\advance\minute by-\hour}
\edef\standardtime{{\ifnum\hour<12 \global\amorpm={am}%
        \else\global\amorpm={pm}\advance\hour by-12 \fi
        \ifnum\hour=0 \hour=12 \fi
        \number\hour:\ifnum\minute<10 0\fi\number\minute\the\amorpm}}
\edef\militarytime{\number\hour:\ifnum\minute<10 0\fi\number\minute}

\def\draftlabel#1{{\@bsphack\if@filesw {\let\thepage\relax
   \xdef\@gtempa{\write\@auxout{\string
      \newlabel{#1}{{\@currentlabel}{\thepage}}}}}\@gtempa
   \if@nobreak \ifvmode\nobreak\fi\fi\fi\@esphack}
        \gdef\@eqnlabel{#1}}
\def\@eqnlabel{}
\def\@vacuum{}
\def\draftmarginnote#1{\marginpar{\raggedright\scriptsize\tt#1}}

\def\draft{\oddsidemargin -.5truein
        \def\@oddfoot{\sl preliminary draft \hfil
        \rm\thepage\hfil\sl\today\quad\militarytime}
        \let\@evenfoot\@oddfoot \overfullrule 3pt
        \let\label=\draftlabel
        \let\marginnote=\draftmarginnote
   \def\@eqnnum{(\theequation)\rlap{\kern\marginparsep\tt\@eqnlabel}%
\global\let\@eqnlabel\@vacuum}  }
\newcounter{app}
\newcounter{sapp}[app]

\newcommand{\app}[1]{
\refstepcounter{app}{\noindent\Large\bf Appendix
 \ #1 \par \vspace{5mm}}
\setcounter{equation}{0}
\def\theequation{\Alph{app}.\arabic{equation}}}
\def\thesapp{\Alph{app}.\arabic{sapp}}
\newcommand{\sapp}[1]{\par \refstepcounter{sapp}{\noindent\large\bf \thesapp
\ #1 \par \vspace{3mm}}
\def\theequation{\Alph{app}.\arabic{equation}}}
\makeatletter
\newdimen\normalarrayskip              % skip between lines
\newdimen\minarrayskip                 % minimal skip between lines
\normalarrayskip\baselineskip
\minarrayskip\jot
\newif\ifold             \oldtrue            
\def\arraymode{\ifold\relax\else\displaystyle\fi} % mode of array entries
\def\eqnumphantom{\phantom{\mbox{\rm
(\theequation)}}}% right phantom in eqnarray
\def\@arrayskip{\ifold\baselineskip\z@\lineskip\z@
     \else
     \baselineskip\minarrayskip\lineskip2\minarrayskip\fi}
\def\@arrayclassz{\ifcase \@lastchclass \@acolampacol \or
\@ampacol \or \or \or \@addamp \or
   \@acolampacol \or \@firstampfalse \@acol \fi
\edef\@preamble{\@preamble
  \ifcase \@chnum
     \hfil$\relax\arraymode\@sharp$\hfil
     \or $\relax\arraymode\@sharp$\hfil
     \or \hfil$\relax\arraymode\@sharp$\fi}}
\def\@array[#1]#2{\setbox\@arstrutbox=\hbox{\vrule
     height\arraystretch \ht\strutbox
     depth\arraystretch \dp\strutbox
     width\z@}\@mkpream{#2}\edef\@preamble{\halign \noexpand\@halignto
\bgroup \tabskip\z@ \@arstrut \@preamble \tabskip\z@ \cr}%
\let\@startpbox\@@startpbox \let\@endpbox\@@endpbox
  \if #1t\vtop \else \if#1b\vbox \else \vcenter \fi\fi
  \bgroup \let\par\relax
  \let\@sharp##\let\protect\relax
  \@arrayskip\@preamble}
\def\eqnarray{\stepcounter{equation}%
              \let\@currentlabel=\theequation
              \global\@eqnswtrue
              \global\@eqcnt\z@
              \tabskip\@centering
              \let\\=\@eqncr
              $$%
 \halign to \displaywidth\bgroup
    \eqnumphantom\@eqnsel\hskip\@centering
    $\displaystyle \tabskip\z@ {##}$%
    &\global\@eqcnt\@ne \hskip 1.2\arraycolsep
         $\displaystyle\arraymode{##}$\hfil
    &\global\@eqcnt\tw@ \hskip 1.2\arraycolsep
         $\displaystyle\tabskip\z@{##}$\hfil
         \tabskip\@centering
    &{##}\tabskip\z@\cr}
\makeatother
\catcode`\@=11
\ifcase\@ptsize
 \font\tenmsa=msam10
 \font\sevenmsa=msam7
 \font\fivemsa=msam5
 \font\tenmsb=msbm10
 \font\sevenmsb=msbm7
 \font\fivemsb=msbm5
 \font\teneu=eufm10
 \font\seveneu=eufm7
 \font\fiveeu=eufm5
 \font\tenib=cmmib10
 \font\sevenib=cmmib7
 \font\fiveib=cmmib5
\or
 \font\tenmsa=msam10 scaled \magstephalf
 \font\sevenmsa=msam7 scaled \magstephalf
 \font\fivemsa=msam5 scaled \magstephalf
 \font\tenmsb=msbm10 scaled \magstephalf
 \font\sevenmsb=msbm7 scaled \magstephalf
 \font\fivemsb=msbm5  scaled \magstephalf
 \font\teneu=eufm10  scaled \magstephalf
 \font\seveneu=eufm7  scaled \magstephalf
 \font\fiveeu=eufm5   scaled \magstephalf
 \font\tenib=cmmib10  scaled \magstephalf
 \font\sevenib=cmmib7  scaled \magstephalf
 \font\fiveib=cmmib5   scaled \magstephalf
\or
 \font\tenmsa=msam10 scaled \magstep1
 \font\sevenmsa=msam7 scaled \magstep1
 \font\fivemsa=msam5  scaled \magstep1
 \font\tenmsb=msbm10 scaled \magstep1
 \font\sevenmsb=msbm7 scaled \magstep1
 \font\fivemsb=msbm5  scaled \magstep1
 \font\teneu=eufm10   scaled \magstep1
 \font\seveneu=eufm7 scaled \magstep1
 \font\fiveeu=eufm5 scaled \magstep1
 \font\tenib=cmmib10     scaled \magstep1
 \font\sevenib=cmmib7   scaled \magstep1
 \font\fiveib=cmmib5   scaled \magstep1
\fi
\newfam\msafam
\textfont\msafam=\tenmsa  \scriptfont\msafam=\sevenmsa
  \scriptscriptfont\msafam=\fivemsa
\newfam\msbfam
\textfont\msbfam=\tenmsb  \scriptfont\msbfam=\sevenmsb
  \scriptscriptfont\msbfam=\fivemsb
\def\Bbb{\ifmmode\let\next\Bbb@\else
 \def\next{\errmessage{Use \string\Bbb\space only in math mode}}\fi\next}
\def\Bbb@#1{{\Bbb@@{#1}}}
\def\Bbb@@#1{\fam\msbfam#1}
\newfam\eufam
\textfont\eufam=\teneu  \scriptfont\eufam=\seveneu
  \scriptscriptfont\eufam=\fiveeu
\def\frak{\ifmmode\let\next\frak@\else
 \def\next{\errmessage{Use \string\frak\space only in math mode}}\fi\next}
\def\frak@#1{{\frak@@{#1}}}
\def\frak@@#1{\fam\eufam#1}
\newfam\ibfam
\textfont\ibfam=\tenib  \scriptfont\ibfam=\sevenib
  \scriptscriptfont\ibfam=\fiveib
\def\bold{\ifmmode\let\next\bold@\else
 \def\next{\errmessage{Use \string\bold\space only in math mode}}\fi\next}
\def\bold@#1{{\bold@@{#1}}}
\def\bold@@#1{\fam\ibfam#1}
\def\hexnumber@#1{\ifcase#1 0\or 1\or 2\or 3\or 4\or 5\or 6\or 7\or 8\or
 9\or A\or B\or C\or D\or E\or F\fi}
\def\newsymbolb#1#2#3#4{\mathchardef#1="#2\hexnumber@\msbfam#3#4}
\def\newsymbola#1#2#3#4{\mathchardef#1="#2\hexnumber@\msafam#3#4}
\relax
\hybrid
\begin{document}
%%%%%%%%%%%%%%%%%%%%%%%%%%%%%%%%%%%%%%%%%%%%%%%%%%%%%%%%%%%%%%%%%%%%%%%%%%
%%%%%%%%%%%%%%%%%  Beginning of the personal definitions  %%%%%%%%%%%%%%%%
%%%%%%%%%%%%%%%%%%%%%%%%%%%%%%%%%%%%%%%%%%%%%%%%%%%%%%%%%%%%%%%%%%%%%%%%%%
\def\bea{\begin{eqnarray}}
\def\eea{\end{eqnarray}}
\def\beq{\begin{equation}}          \def\bn{\beq}
\def\eeq{\end{equation}}            \def\ed{\eeq}
\def\nn{\nonumber}                  \def\g{\gamma}
\def\Uq{U_q(\widehat{\frak{sl}}_2)}
\def\Uqp{U_q(\widehat{\frak{sl}}'_2)}
\def\Uqd{U^{*}_q(\widehat{\frak{sl}}_2)}
\def\uq{U_q({sl}_2)}
\def\uqd{U^*_q({sl}_2)}
\def\slaff{\frak{sl}^\prime_2}
\def\aff{\widehat{\frak{sl}}_2}
\def\ot{\otimes}
\def\id{\mbox{\rm id}}
\def\Re{{\rm Re}\,}
\def\RR{\Bbb{R}}               %%%%%   This four lines should be
\def\ZZ{\Bbb{Z}}              %%%%%         corrected to the AMS fonts
\def\CC{\Bbb{C}}                %%%%%%%
\def\r#1{\mbox{(}\ref{#1}\mbox{)}}
\def\d{\delta}
\def\D{\Delta}
\def\da{{\partial_\alpha}}
\let\da=p
\def\R{{\cal R}}
\def\h{\hbar}
\def\Ga#1{\Gamma\left(#1\right)}
\def\ep{\varepsilon}
\def\ve{\ep}
\def\fract#1#2{{\mbox{\footnotesize $#1$}\over\mbox{\footnotesize $#2$}}}
\def\stackreb#1#2{\ \mathrel{\mathop{#1}\limits_{#2}}}
\def\res#1{\stackreb{\mbox{\rm res}}{#1}}
\def\Res#1{\stackreb{\mbox{\rm Res}}{#1}}
\let\dis=\displaystyle
\def\ee{{\rm e}}
\def\D{\Delta}
\renewcommand{\theequation}{\mbox{\rm {\thesection}.{\arabic{equation}}}}
\newsymbolb{\ltimes}26E
\newsymbolb{\rtimes}26F
\newsymbola{\probaF}23F
\def\Y-{\widehat{Y}^-}
\font\fraksect=eufm10 scaled 1440
\def\DYsect{\widehat{DY(\hbox{\fraksect sl}_2)}}
\def\DY{\widehat{DY(\frak{sl}_2)}}
\def\Yd{\DY}
\def\Ydd{\DY}
\let\z=z
\let\b=z
\let\u=u
\let\v=v
\let\w=\zeta
%%%%%%%%%%%%%%%%%%%%%%%%%%%%%%%%%%%%%%%%%%%%%%%%%%%%%%%%%%%%%%%%%%%%%%%%%%
%%%%%%%%%%%%%%%%%%%%%%%  End of the personal definitions  %%%%%%%%%%%%%%%%
%%%%%%%%%%%%%%%%%%%%%%%%%%%%%%%%%%%%%%%%%%%%%%%%%%%%%%%%%%%%%%%%%%%%%%%%%%
\begin{titlepage}
\begin{center}
\hfill DFTUZ/95/28\\
\hfill ITEP-TH-15/95\\
%\hfill JINR-95-??\\
\hfill q-alg/9602030\\
\bigskip\bigskip
{\Large\bf  Intertwining Operators for the Central Extension\\
of the Yangian Double}\\
\bigskip
\bigskip
{\large S. Khoroshkin\footnote{E-mail: khoroshkin@vitep1.itep.ru},
D. Lebedev\footnote{E-mail: lebedev@vitep1.itep.ru}}\\
\medskip
{\it Institute of Theoretical \& Experimental Physics\\
117259 Moscow, Russia}\\
\bigskip
{\large S. Pakuliak}\footnote{E-mail: pakuliak@thsun1.jinr.dubna.su}\\
\bigskip
{\it Departamento de F\'\i sica Te\'orica, Facultad de Ciencias\\
Universidad de Zaragoza, 50009 Zaragoza, Spain}\\
\medskip
{\it and}\\
\medskip
{\it Bogoliubov Laboratory of Theoretical Physics, JINR\\
141980 Dubna, Moscow region, Russia}\\
\bigskip
\bigskip
\bigskip
%{Revised \today}
\end{center}
\begin{abstract}
We continue the investigation of the central extended Yangian double \cite{K}.
In this paper we study the intertwining operators for  certain
infinite dimensional representations of $\Yd$, which are deformed analogs
of the highest weight representations of the  affine algebra
$\aff$ at level 1.
We give bosonized expressions for intertwining operators, verify that they
generate an algebra isomorphic to Zamolodchikov--Faddeev algebra
for the $SU(2)$-invariant Thirring model. We compose from them
$L$-operators by Miki's prescription and verify that they coincide
with $L$-operators constructed from universal ${\cal R}$-matrix.
\end{abstract}
\end{titlepage}
\clearpage
\newpage

\setcounter{equation}{0}
\setcounter{footnote}{0}

\section{Introduction}

The Yangian $Y(g)$ was introduced by V.G. Drinfeld \cite{D1}. As a
Hopf algebra it is a deformation of the universal enveloping algebra
 $U(g[u])$ of $g$-valued
polynomial currents where $g$ is a simple Lie algebra.
Recently, it was understood that in physical applications one needs
 the double of the Yangian \cite{S2}. In
\cite{KT} the double of the Yangian was described quite explicitely in
 terms of  Drinfeld's construction of a
quantum double  \cite{D3}.  After adding to the Yangian a
derivation operator, the same quantum double construction produces a central
extension of the Yangian double \cite{K}. The  latter algebra has a richer
representation theory than just the
 Yangian double has. For example,  infinite
dimensional representations which are deformed analogs of the highest weight
representations of the affine algebra can be constructed. This type of
representations appears to be very important in understanding the structure
of conformal field theories associated with affine Lie algebras and the
structure of the quantum integrable models associated with deformation of the
universal enveloping affine algebras.  In the  paper \cite{XXZ}
and then in the book \cite{JM} the representation theory of the infinite
dimensional representations of $\Uq$ has been used in order to solve
completely the quantum spin $1/2$ lattice XXZ model.
A complete solution of the
quantum model means the possibility to write down explicitly an analytical
integral representation for arbitrary form-factor of any local operator
in the model.

Some years before  an  approach to relativistic integrable
massive field theories have been developed in the works by F.A. Smirnov
\cite{S1}. This approach did not use representation theory of the
hidden non-abelian symmetry of the model which ensures  its integrability,
but the ideas of phenomenological bootstrap. For some quantum integrable
models such as Sin-Gordon, SU(2)-invariant Thirring, etc. models
the bootstrap program have led to a  complete solution.
Understanding the fact
 that the dynamical non-abelian symmetry algebra for
SU(2)-invariant Thirring model is the Yangian double (see \cite{S2} and
reference therein) opens a possibility to describe the structure of
the model (local operators, Zamolodchikov--Faddeev operators which
creates particles, etc.) in terms of the representation theory of the
Yangian double. It appears that it should be central extended
Yangian double.

On the other hand, a free field approach to massive integrable field
theories has been developed in \cite{L}. The starting point
there has been the Zamolodchikov--Faddeev (ZF) algebra \cite{ZZ,F} for the
operators which create and annihilate particles in the model
and describe the local operators.   The
author of this paper was able to construct a free field
representation for ZF operators
using an ultraviolet
regularization of the Fock modules.
Our approach  shows that there
exists possibility to construct a free field representation of
the ZF operators in massive integrable field
theories where only a minor regularization in the definition of
these operators will be necessary
(see \r{Phi-}--\r{comparII}).

The paper is organized as follows. After the definition of $\Yd$ in section 2,
an infinite dimensional representations of this algebra are constructed
in terms of free field in section 3. The next two sections are
 devoted to the free field representation
 of the intertwining operators and to their
commutation and normalization relations.
In section 6  an universal ${\cal R}$-matrix description
of the central extended Yangian double
and Ding--Frenkel isomorphism \cite{DF} between the ``new''
\cite{D2} and ``$RLL$'' formulations \cite{FRT,RS}
of $\Yd$  is given \cite{KT,K}.
In the last section  we consider
 a free field construction of the
$L$-operators corresponding
to the  central extension of the Yangian double at level 1
(Miki's formulas \cite{M}) and verify that they coincide with $L$-operators
constructed from universal ${\cal R}$-matrix.
Some explicit calculations are gathered in  Appendix.

%\newpage

\setcounter{equation}{0}

\section{Central Extension of $\DYsect$}

\let\hsp=\qquad
The Yangian $Y(sl_2)$ is a Hopf algebra generated by the
elements $e_k$, $f_k$, $h_k$, $k\geq0$, subject to the relations
\bea
&[h_{k},h_{l}]\ = \ 0 \ , \qquad  [e_{k}, f_{l}] \ =\ h_{k+l}\ ,\nn\\
&[h_{0},e_{l}]\ = \ 2e_{l}\ , \qquad [h_{0},f_{l}] = -2f_{l}\ , \nn\\
&[h_{k+1},e_{l}]\ -\ [h_{k},e_{l+1}]\ = \
\h\{ h_{k},e_{l}\}\ ,\nn\\
&[h_{k+1},f_{l}]\ -\ [h_{k},f_{l+1}]\ = \
-\h\{h_{k},f_{l}\}\ ,\nn\\
&[e_{k+1},e_{l}] \ -\ [e_{k},e_{l+1}]\ = \
\h\{ e_{k},e_{l}\}\ ,\nn\\
&[f_{k+1},f_{l}] \ - \ [f_{k},f_{l+1}]\ =      \
-\h\{f_{k},f_{l}\}\ ,
\label{2.1}
\eea
 where $\h$ is a deformation parameter  and $\{ a,b\} =ab+ba$.
The coalgebra structure is uniquely defined by the relations
 $$\D(e_0)= e_0\otimes 1+1\otimes e_0\ , \hsp
 \D(h_0)= h_0\otimes 1+1\otimes h_0\ , \hsp
 \D(f_0)= f_0\otimes 1+1\otimes f_0\ ,$$
 $$\D(e_1)= e_1\otimes 1+1\otimes e_1 +\h h_0\otimes e_0\ ,\hsp
 \D(f_1)= f_1\otimes 1+1\otimes f_1 +\h f_0\otimes h_0\ ,$$
 $$
 \D(h_1)= h_1\otimes 1+1\otimes h_1+\h h_0\otimes h_0-2f_0\otimes
 e_0\ .
 %\label{2.2}
 $$
In the form \r{2.1}
  the Yangian $Y(sl_2)$
 appeared in \cite{D2} in so called ``new'' realization
which is a deformed analog of the loop realization of the affine algebras.
 In this paper we study representations and intertwining operators for
the central extension $\DY$ of the
quantum double of $Y(sl_2)$. The Hopf algebra
 $\DY$ was introduced in \cite{K}. It can be described as a formal
 completion  of  algebra with generators
  $d$, central
 element $c$ and $e_k, f_k, h_k$, $k\in {\ZZ}$,
gathered into generating functions
 \beq
\label{2.0}
e^{\pm}(\u)=\pm\sum_{k\geq0\atop k<0} e_k \u^{-k-1},
\hsp
f^{\pm}(\u)=\pm\sum_{k\geq0\atop k<0} f_k \u^{-k-1},
\hsp
h^{\pm}(\u)=1\pm\h\sum_{k\geq0\atop k<0}h_{k}\u^{-k-1},
\eeq
 $$ e(\u)= e^+(\u)-e^-(\u),\hsp  f(\u)= f^+(\u)-f^-(\u),$$
with the relations
$$[d,e(\u)]=\frac{\mbox{d}}{\mbox{d}u}e(\u), \hsp
 [d,f(\u)]=\frac{\mbox{d}}{\mbox{d}u}f(\u), \hsp
 [d,h^\pm(\u)]=\frac{\mbox{d}}{\mbox{d}\u}h^\pm(\u),$$
\begin{eqnarray}
e(\u)e(\v)&=&{\u-\v+\h\over \u-\v-\h}\ e(\v)e(\u) \nn\\
f(\u)f(\v)&=&{\u-\v-\h\over \u-\v+\h}\ f(\v)f(\u) \nn\\
h^\pm(\u)e(\v)&=&{\u-\v+\h\over \u-\v-\h}\ e(\v)h^\pm(\u) \nn\\
h^+(\u)f(\v)&=&{\u-\v-\h-\h c\over \u-\v+\h-\h c}\ f(\v)h^+(\u) \nn\\
h^-(\u)f(\v)&=&{\u-\v-\h\over \u-\v+\h}\ f(\v)h^-(\u) \nn\\
h^+(\u)h^-(\v)
&=&{\u-\v+\h\over \u-\v-\h}\cdot{\u-\v-\h-\h c\over \u-\v+\h-\h c}\
h^-(\v)h^+(\u)
\nn\\
{[}e(\u),f(\v){]}&=&
{1\over \h}\left(\delta(\u-(\v+\h c))h^+(\u)-\delta(\u-\v)h^-(\v) \right)
\label{DY2}
\end{eqnarray}
where
$$\delta(\u-\v)=\sum_{n+m=-1}\u^n\v^m, \hsp
\delta(\u-\v)g(\u)=\delta(\u-\v)g(\v).$$

 This algebra admits a filtration
\beq
\ldots\subset C_{-n}\subset \ldots\subset C_{-1} \subset C_{0}
\subset C_{1}\ldots\subset C_{n}\ldots \subset C
\label{2.7}
\eeq
defined by the conditions $\mbox{deg}\, e_{k}=\mbox{deg}\, f_{k}=
\mbox{deg}\, h_{k}=k$; $\deg \left\{ x\in C_m\right\}\leq m$.
Then $\DY$ is a formal completion of $C$ with respect to filtration
\r{2.7}.

The comultiplication in $\DY$ is given by the relations
$$
\D(h^{\varepsilon}(\u))\ =\
\sum_{k=0}^\infty (-1)^k(k+1)\h^{2k}
\big(f^{\ve}(\u+\h-\d_{\varepsilon ,+}\h c_1)\big)^kh^{\ve}(\u)\otimes
h^{\ve}(\u-\d_{\varepsilon ,+}\h c_1) \big(e^{\ve}(\u+\h-
\d_{\varepsilon ,+}\h c_1)\big)^{k},
$$
\bea
\D(e^{\varepsilon}(\u))&=&e^{\ve}(\u)\otimes 1
+\sum_{k=0}^\infty (-1)^k\h^{2k}
\big(f^{\ve}(\u+\h-\d_{\varepsilon ,+}\h c_1)\big)^k
h^{\ve}(\u)\otimes\big(e^{\ve}(\u-\d_{\varepsilon ,+}\h c_1)
\big)^{k+1},\nn      \\
\D(f^{\varepsilon}(\u))&=&1\otimes f^{\ve}(\u)
+\sum_{k=0}^\infty (-1)^k\h^{2k}
\big(f^{\ve}(\u+\d_{\varepsilon ,+}\h c_2)\big)^{k+1}\otimes
h^{\ve}(\u)\big(e^{\ve}(\u+\h)\big)^k,
        \label{3.6}
        \eea
        where  $\varepsilon =\pm,\; \d_{+ ,+}=1$, $\d_{- ,+}=0$
and $c_1$, $c_2$ means $c\ot 1$, $1\ot c$ respectively.

The Hopf algebra $\DY$ was derived in two steps. First, it was proved
 in \cite{KT} that the quantum double $\Ydd$ of the Yangian
$Y(\frak{sl}_2)$ can be
 described by formulas \r{DY2}-\r{3.6} with $c=0$ and $d$ being
 dropped.
  The double $\Ydd$ containes Yangian $Y(\frak{sl}_2)$
and its dual with the opposite
  comultiplication as Hopf subalgebras.
  Subalgebra $Y^+=Y(\frak{sl}_2)\subset DY(\frak{sl}_2)$ is generated by
the components of $e^+(\u)$, $f^+(\u)$, $h^+(\u)$ and its dual with
the opposite
comultiplication $Y^-=(Y(\frak{sl}_2))^0$ is a formal completion (see
\r{2.7}) of the subalgebra generated by the components of $e^-(\u)$,
 $f^-(\u)$, $h^-(\u)$.
 The Hopf pairing $<,>:\; Y^+\otimes (Y^-)^0\rightarrow \CC$  has the
following description.

 Let $E^\pm$, $H^\pm$, $F^\pm$ be subalgebras  (or
 their completions in $Y^-$ case) generated by the components of
  $e^\pm(u)$, $h^\pm(u)$, $f^\pm(u)$. Subalgebras $E^\pm$ and
 $F^\pm$ do not contain the unit.
 Then
 \bn
 {<}e^+h^+f^+,f^-h^-e^-{>}= {<}e^+,f^-{>}{<}h^+,h^-{>}{<}f^+,e^-{>}
 \label{factorization}
 \ed
 for any elements $e^\pm\in E^\pm$, $h^\pm\in H^\pm$,
 $f^\pm\in F^\pm$
 and
 $${<}e^+(\u),f^-(\v){>}= \frac{1}{\h (\u-\v)}\ ,\hsp
 {<}f^+(\u),e^-(\v){>}= \frac{1}{\h (\u-\v)}\ ,$$
 $${<}h^+(\u),h^-(\v){>}=\frac{\u-\v+\h}{\u-\v-\h}\ .$$
 These properties defines the pairing uniquely.

 At the second step the Hopf algebra $Y^-$ is extended to a semidirect
 $\Y-=Y^-\rtimes{\CC}[[d]]$, where $d$
  is a derivation
\bn
[d,e^\pm(\u)]=\frac{{\rm d}}{{\rm d}\u}e^\pm(\u)\ , \hsp
[d,h^\pm(\u)]=\frac{{\rm d}}{{\rm d}\u}h^\pm(\u)\ , \hsp
[d,f^\pm(\u)]=\frac{{\rm d}}{{\rm d}\u}f^\pm(\u)\ .
\label{3.0}
\ed
 $\Y-$ is again a Hopf algebra if we put
$$\D(d)=d\otimes 1 +1\otimes d\ .$$
  Finally,  the central extension $\DY$ is the quantum double of $\Y-$.
 Central element $c$ is dual to $d$: $\; {<}c,d{>}=\frac{1}{\h}.$
 The relations \r{DY2} result from the
  general commutation relation in Drinfeld's construction of the
 quantum double \cite{D3}.
Together with  \r{3.6} they define the Hopf algebra
 $\DY$ in the ``new'' realization.
 Concrete calculations are exhibited in
 \cite{K}.

Before passing to the next section where we will give a free field
realization of the commutation relations \r{DY2} at the level $c=1$
we would like to add two remarks.

 One can note first that $\DY$ is a
 deformation of the enveloping algebra $U(\frak{sl}_2[t,t^{-1}]]+{\CC}c+{\CC}
\frac{d}{dt})$ (which differs from $U(\widehat{\frak{sl}}_2)$
by grading element)
in the direction of
rational $r$-matrix $r(u)=\frac{g_a\otimes g^a}{u}$. The
 deformation of Cartan fields are generating
functions $\kappa^\pm(u)={1\over\h}\ln\, h^\pm(u)$ (see \r{Fren-Kac-boson}).

\let\bet=b
 The second remark is about the non-symmetry
of the relations \r{DY2}. One can easily note that the replacement
\bea
h^+(u)&\to&\tilde h^+(u)=h^+(u+(\bet+1/2)c\h)\ ,\nn\\
h^-(u)&\to&\tilde h^-(u)=h^-(u+\bet c\h)\ ,\nn\\
e(u)&\to&\tilde e(u)=e(u+(\bet+1/4)c\h)\ ,\nn\\
f(u)&\to&\tilde f(u)=f(u+(\bet-1/4)c\h)\ ,
\nn%\label{replacement}
\eea
where $\bet$ is  an arbitrary complex number,
will move the commutation relations \r{DY2} into more symmetric
form
\begin{eqnarray}
\tilde
h^+(\u)\tilde h^-(\v)&=&{\u-\v+\h+\h c/2\over \u-\v-\h-\h c/2}\cdot
{\u-\v-\h-\h c/2\over \u-\v+\h-\h c/2}\
\tilde
h^-(\v)\tilde h^+(\u)\ ,\nn\\
\tilde h^\pm(\u)\tilde e(\v)&=&{\u-\v+\h\pm\h c/4\over \u-\v-\h\pm\h c/4}\
\tilde e(\v)\tilde h^\pm(\u)\ , \nn\\
\tilde h^\pm(\u)\tilde f(\v)&=&{\u-\v-\h\mp\h c/4\over \u-\v+\h\mp\h c/4}\
\tilde f(\v)\tilde h^\pm(\u)\ , \nn\\
\tilde e(\u)\tilde f(\v)-\tilde f(\v)\tilde e(\u) &=&
{1\over \h}\left(\delta\left(\u-\v-{\h c\over 2}\right)
\tilde h^+\left(\u-{\h c\over 4}\right)\right.\nn\\
&-&\left.
\delta\left(\u-\v+{\h c\over2}\right)\tilde
h^-\left(\v-{\h c\over 4}\right) \right)\ .
\nn%\label{comDYsymm}
\end{eqnarray}
For quantum affine algebras such a renormalization have been initiated by the
condition of consistency with Cartan involution. There is no such motivation
in our case so we do not use this replacement any more.

\setcounter{equation}{0}

\section{Bosonization Formulas for Level 1}

Let ${\cal H}$ be the Heisenberg algebra generated by free
bosons with zero modes $a_{\pm n}$, $n=1,2,\ldots$,
$a_0$, $\da$  with commutation relations
$$%\begin{equation}
%\label{bosons}
[a_n,a_m] =  n \delta_{n+m,0}\ ,\quad [\da,a_0]=2\ .
$$%\end{equation}
Let $V_i, i=0,1$ be the formal power series extensions of the Fock spaces
\beq
V_i={\CC}[[a_{-1},\ldots ,a_{-n},\ldots ]]\otimes
        \left(\oplus_{n\in{\ZZ}+i/2}{\CC}\ee^{na_0}\right)
\label{modules}
\eeq
with the action of bosons on these spaces
\bea
a_n&=&\hbox{the left multiplication by $a_n\otimes1$\ \  for $n<0$}\ ,\nn\\
   &=&{[}a_n,\ \cdot\ {]}\otimes 1\quad \hbox{for $n>0$}\ ,\nn\\
\ee^{n_1a_0}(a_{-j_k}\cdots a_{-j_1}\otimes\ee^{n_2a_0})&=&
a_{-j_k}\cdots a_{-j_1}\otimes\ee^{(n_1+n_2)a_0}\ ,\nn\\
u^{\da} (a_{-j_k}\cdots a_{-j_1}\otimes\ee^{na_0})&=&u^{2n}
a_{-j_k}\cdots a_{-j_1}\otimes\ee^{na_0}\ .\nn%\label{action}
\eea
The End$V_i$-valued generating functions (fields) \cite{K}
\begin{eqnarray}
e(\u)&=&
\exp\left(\sum_{n=1}^{\infty} {a_{-n}\over n}\left[(\u-\h)^n+\u^n\right]
\right)
\ee^{a_0}\u^{\da}
\exp\left(-\sum_{n=1}^{\infty} {a_{n}\over n}\u^{-n} \right)
\ , \nn\\
f(\u)&=&
\exp\left(-\sum_{n=1}^{\infty} {a_{-n}\over n}\left[(\u+\h)^n+\u^n\right]
\right)
\ee^{-a_0}\u^{-\da}
\exp\left(\sum_{n=1}^{\infty} {a_{n}\over n}\u^{-n} \right)
\ , \nn\\
h^-(\u)&=&
\exp\left(\sum_{n=1}^{\infty} {a_{-n}\over n}\left[(\u-\h)^n-(\u+\h)^n\right]
\right)\ ,\nn\\
h^+(u)&=&
\exp\left(\sum_{n=1}^{\infty} {a_{n}\over n}\left[(\u-\h)^{-n}-\u^{-n}\right]
\right)
\left({\u-\h\over \u}\right)^{-\da}
                                \label{bosonization}
\end{eqnarray}
satisfy commutation relations \r{DY2} with $c=1$.

        In terms of generating  functions
        \bea
        a_+(z)&=&\sum_{n\geq 1}\frac{a_n}{n}z^{-n}-p\log z\ ,\hsp
        a_-(z)=\sum_{n\geq 1}\frac{a_{-n}}{n}z^{n}+\frac{a_0}{2}\ ,
        \label{5.2} \\
        a(z)&=&a_+(z)-a_-(z),\hsp \phi_\pm(z)=\exp a_\pm(z)\ ,
        \label{5.2a}    \\
        {[}a_+(z),a_-(y){]}&=&-\log(z-y), \hsp |y|<|z|\ ,
        \label{5.2b}
        \eea
 a free filed realization \r{bosonization} looks as follows:
\bea
e(\u)&=&\phi_-(\u-\h)\phi_-(\u)\phi_+^{-1}(\u)\ ,\hsp
 f(\u)=\phi_-^{-1}(\u+\h)\phi_-^{-1}(\u)\phi_+(\u)\ ,\nn\\
h^+(\u)&=&\phi_+(\u-\h)\phi_+^{-1}(\u)\ ,\hsp
 h^-(\u)=\phi_-(\u-\h)\phi_-^{-1}(\u+\h)\ ,
 \label{h-fields}                     \\
 \ee^{\g d}\phi_\pm(\u)&=&\phi_\pm(\u+\g)\ee^{\g d}\ ,
\hsp \ee^{\g d}(1\otimes 1)
 =1\otimes 1. \nn
\eea

One can note that in the classical limit ($\hbar\to 0$) the bosonization
formulas goes to the Frenkel-Kac (homogeneous) realization of the affine
algebra $\widehat{\frak{sl}}_2$ at level 1
\cite{FK}
\begin{eqnarray}
e(\u)&=&
\exp\left(\sum_{n=1}^{\infty} {2a_{-n}\over n}\u^n\right)
\ee^{a_0}\u^{\da}
\exp\left(-\sum_{n=1}^{\infty} {a_{n}\over n}\u^{-n} \right)
\ , \nn\\
f(\u)&=&
\exp\left(-\sum_{n=1}^{\infty} {2a_{-n}\over n}\u^n\right)
\ee^{-a_0}\u^{-\da}
\exp\left(\sum_{n=1}^{\infty} {a_{n}\over n}\u^{-n} \right)
\ , \nn\\
\kappa(\u)&=&\kappa^+(\u)-\kappa^-(\u)=
\sum_{n=1}^{\infty} 2a_{-n} \u^{n-1} +
\sum_{n=1}^{\infty} a_{n} \u^{-n-1} + \da \u^{-1}\ .
\label{Fren-Kac-boson}
\end{eqnarray}

Let us discuss now the structure of the representation spaces
$V_i$, $i=0,1$. In the case of affine algebra
$\widehat{\frak{sl}}_2$ or quantum
(trigonometric)
deformation of the affine algebra $\Uq$ there is a finite number
of Chevalley generators which generate all the infinite dimensional
algebra by means of the multiple commutators.
It appears also that Chevalley generators coincide with
certain components of the currents that form the loop (the new)
realization of
affine (quantum affine) algebra. This gives a combinatorial description
of the representations
$V_i$, $i=0,1$ with highest weight vectors
$u_0=1\ot 1$, $u_1=1\ot \ee^{a_0/2}$ in terms of successive
 finite-dimensional modules of the two embedded copies of
$U_{q}(\frak{sl}_2)$.

 In the case of $\DY$ the second copy of $\frak{sl}_2$
is missing, and, on the other hand,
the action of the components of the currents $e_k$ or $f_k$ onto any element
of the representation spaces  $V_i$  produces infinite sum of the
elements of the same sort. For example,
$$%\beq
e_{-1} 1\ot 1 = \exp\left(\sum_{n=1}^\infty {a_n \hbar^n\over n}\right)
\ot \ee^{a_0} \ .
%\label{example}
$$%\eeq
It means that in case of the central extension of the Yangian double
 we have no analogous combinatorial description of the representation spaces
$V_0$ and $V_1$ .  Nevertheless these spaces
are properly defined together with an action  of $\DY$ on them.

One can show that there exists an outer automorphism $\nu$ of $\Yd$
such that $V_i$ is isomorphic to $V_{1-i}^\nu$, $i=0,1$.

To construct intertwining operators we  need also
the finite-dimensional evaluation modules for $\Yd$.
 The following
formulas can be found in \cite{KT}
\begin{eqnarray}
h^+(\u)v_{\pm,\b}&=&{\u-\b\pm \h\over \u-\b}
 v_{\pm,\b},\quad |\b|<|\u|\ ,\nn\\
h^-(\u)v_{\pm,\b}&=&{\u-\b\pm \h\over \u-\b}
 v_{\pm,\b},\quad |\b|>|\u|\ ,\nn\\
e^\pm(\u)v_{+,\b}&=&0,\quad f^\pm(\u)v_{-,\b}\, =\, 0\ ,       \nn\\
e^+(\u)v_{-,\b}&=&{1\over \u-\b}
 v_{+,\b},\quad
f^+(\u)v_{+,\b}\, =\, {1\over \u-\b}
 v_{-,\b},\quad
|\b|<|\u|\ ,\nn\\
e^-(\u)v_{-,\b}&=&{1\over \u-\b}
 v_{+,\b},\quad
f^-(\u)v_{+,\b}\, =\, {1\over \u-\b}
 v_{-,\b},\quad
|\b|>|\u|\ ,
\label{eval}
\end{eqnarray}
where $v_{\pm,\b}$ denotes the elements of the
evaluation module
$V_\b=\{V\otimes \CC[\b,\b^{-1}]\}$ and $V=\{v_+,v_-\}$.

\setcounter{equation}{0}
\section{Intertwining Operators}

Let us define the following intertwining
operators
\begin{eqnarray}
&\Phi^{(i)}(\z): V_i\to V_{1-i}\otimes V_\z,\quad
\Phi^{*(i)}(\z):V_i\otimes V_\z\to V_{1-i}
\ ,\label{typeI}\\
&\Psi^{(i)}(\b): V_i\to V_\b\otimes V_{1-i},\quad
\Psi^{*(i)}(\b):V_\b\otimes V_i\to V_{1-i}
\label{typeII}
\end{eqnarray}
which commute with the action of the Yangian double
\begin{eqnarray}
\Phi^{(i)}(\z) x &=& \Delta(x) \Phi^{(i)}(\z),\quad
\Phi^{*(i)}(\z)\Delta(x)= x\Phi^{*(i)}(\z),\nn\\
\Psi^{(i)}(\b) x &=& \Delta(x) \Psi^{(i)}(\b),\quad
\Psi^{*(i)}(\b)\Delta(x)= x\Psi^{*(i)}(\b),\quad
\forall x\in\Yd \ .
                          \label{inteq}
\end{eqnarray}
The components of the intertwinning operators are defined as follows
\bea%$$%\begin{equation}
\Phi^{(i)}
(\z) v  &=&
\Phi^{(i)}_+(\z)v\otimes v_+ + \Phi^{(i)}_-(\z)v\otimes v_-\ ,\quad
\Phi^{*(i)}(\z)(v\otimes v_\pm)= \Phi^{*(i)}_\pm(\z)v \ ,
  \nn\\% \label{components}
%$$%\end{equation}
%$$%\begin{equation}
\Psi^{(i)}(\b) v  &=&
v_+\otimes\Psi^{(i)}_+(\b)v + v_-\otimes\Psi^{(i)}_-(\b)v\ ,\quad
\Psi^{*(i)}(\b)(v_\pm\otimes v)= \Psi^{*(i)}_\pm(\b)v \ ,
 \nn  %\label{componentsII}
\eea%$$%\end{equation}
where $v\in V_i$.

According to the terminology proposed in \cite{JM} we call the
intertwining operators \r{typeI} of type I operators and the
operators \r{typeII} type II ones.
Let us also fix
 the normalization of the intertwining operators which yields
the dependence of these operators on the index $i$
\bea\label{normaliz}
\Phi^{(i)}(\z)u_{i}&=&(-\z)^{-i/2}
\ep_{1-i}u_{{1-i}}
\otimes v_{\ep_i}+\cdots,\nn\\
\Psi^{*(i)}(\b)(v_{\ep_{1-i}}\ot u_{i})&=&(-\b)^{-i/2}u_{{1-i}}+\cdots,\quad
\ep_0=-,\ \ \ep_1=+\  ,
\eea
where we denote
$$u_0=1\ot 1\qquad \hbox{and} \qquad u_1=1\ot \ee^{a_0/2} $$
and dots in \r{normaliz} mean the terms containing positive
powers of $\hbar$. The normalization \r{normaliz}
is not standard\footnote{In contrast to the normalization
of the   intertwining operators
for  the quantum affine algebra $\Uq$,
the r.h.s. of \r{normaliz} contains both
 positive and negative powers of the spectral
parameter. The terms with negative  powers of the spectral
parameter disappear in the
classical limit $\h\to 0$.}
but it allows us to write down precise expressions for
$$
\Phi_\ep=\Phi^{(0)}_\ep\oplus\Phi^{(1)}_\ep:V_0\oplus V_1\rightarrow
V_1\oplus V_0
\quad\hbox{and}\quad
\Psi^*_\ep=\Psi^{*(0)}_\ep\oplus\Psi^{*(1)}_\ep:V_0\oplus V_1\rightarrow
V_1\oplus V_0
$$
without dependence on the index $i$.

Let $\eta_+(z)$ be the following $\mbox{End}\, V_i$-valued function
\beq
\eta_+(z)=
\stackreb{\mbox{\rm lim}}{K\to \infty}
(2\h K)^{-\da/2}
\prod_{k=0}^K {\phi_+(z-2k\h)\over\phi_+(z-\h-2k\h)}\ .
\label{eta}
\eeq
We have the following
\medskip

\noindent
{\bf Proposition 1.} {\it Intertwining operators \r{typeI}
 have the free field realiza\-ti\-on:}
\begin{eqnarray}
\Phi_-(\z)&=&
\phi_-(\z+\h)\eta_+^{-1}(\z)\ ,\label{Phi-}\\
\Phi_+(\z)&=&
\Phi_-(\z)f_0-f_0\Phi_-(\z)
=-\h\int_{C}{d\u\over 2i\pi}
{{:}\Phi_-(\z)f(u){:}\over(\u-\z)(\u-\z-\h)}\ ,
\label{Phi+}
\\
   \Phi^{*(i)}_\ep(\z)&=&\ep(-1)^i\Phi^{(i)}_{-\ep}(\z-\h)\ ,
\label{comparI}\\
\Psi^{*}_-(\b)&=&
\phi_-^{-1}(\b)\eta_+(\b)\ ,\label{Psi*-}\\
\Psi^*_+(\b)
&=& e_0 \Psi^*_-(\b)-\Psi^*_-(\b) e_0
=-\h\int_{\tilde C}{d\v\over 2i\pi}
{{:}\Psi^*_-(\z)e(v){:}\over(\v-\b)(\v-\b-\h)}\ ,
\label{Psi+}
\\
   \Psi^{(i)}_\ep(\b)&=&\ep(-1)^{-i}\Psi^{*(i)}_{-\ep}(\b-\h) ,\quad
   \ep=\pm\ , \label{comparII}
\end{eqnarray}
{\it where contours $C$ and $\tilde C$ are such that
the points $\z+\h,0$ ($\b,0$) should be inside the contour $C$ ($\tilde C$)
while the point $\z,\infty$ ($\b+\h,\infty$)
should be  outside  $C$ ($\tilde C$)
respectively.}
\smallskip

To prove this proposition we
need the following
partial knowledge of the
coalgebra structure of the central extended Yangian double $\Yd$
(see also \r{3.6})
\begin{eqnarray}
\Delta h^+(u)&=& h^+(u)\otimes h^+(u-\h c_1)+ F\otimes E \ ,\nn\\
\Delta h^-(u)&=& h^-(u)\otimes h^-(u)+F\otimes E\ ,\nn\\
\Delta e(u)&=&e(u)\otimes 1+ F'\otimes E\ ,\nn\\
\Delta f(u)&=&1\otimes f(u)+ F\otimes E'\ ,\label{comul1}
\end{eqnarray}
and the action of $\Yd$ onto evaluation module \r{eval}.
In \r{comul1} $E$, $F$, $E'$ and $F'$ are generated by the elements
$\{e_k\}$, $\{f_k\}$, $\{e_k,h^\pm_k\}$, $\{f_k,h^\pm_k\}$
respectively.

Taking into account
\r{comul1}
we can obtain   from \r{inteq} for $x=h^+(u)$, $h^-(u)$,
$e(u)$ and $f_0$  the following equations for $\Phi_-^{(i)}(z)$
\begin{eqnarray}
\Phi^{(i)}_-(\z)h^+(u)&=&{u-\z-2\h\over u-\z-\h}h^+(u)\Phi^{(i)}_-(\z)\ ,
\label{fI1}\\
\Phi^{(i)}_-(\z)h^-(u)&=&{u-\z-\h\over u-\z}h^-(u)\Phi^{(i)}_-(\z)\ ,
\label{fI2}\\
\Phi^{(i)}_-(\z)e(u)&=&e(u)\Phi^{(i)}_-(\z)\ ,          \label{fI3}\\
\Phi^{(i)}_-(\z)&=&\Phi^{(i)}_+(\z)f_0- f_0\Phi^{(i)}_+(\z)\ .
           \label{4.16}
\end{eqnarray}
Analogous formulas for the dual operators $\Phi^*(\z)$ are
\begin{eqnarray}
\Phi^{*(i)}_+(\z)h^+(u)&=&{u-\z-\h\over u-\z}h^+(u)\Phi^{*(i)}_+(\z)\ ,
\nn%\label{fI*1}
\\
\Phi^{*(i)}_+(\z)h^-(u)&=&{u-\z\over u-\z+\h}h^-(u)\Phi^{*(i)}_+(\z)\ ,
\nn%\label{fI*2}
\\
\Phi^{*(i)}_+(\z)e(u)&=&e(u)\Phi^{*(i)}_+(\z)\ ,
               \nn%    \label{fI*3}
\\
\Phi^{*(i)}_-(\z)&=&f_0\Phi^{*(i)}_+(\z)- \Phi^{*(i)}_+(\z)f_0\ .
         \nn%    \label{fI*4}
\end{eqnarray}
Comparing  the above formulas we obtain \r{comparI}.
With this choice of the $z$-independent factor $(-1)^i$ in \r{comparI}
we will have  the same normalization
condition for the dual intertwining
 operators as for the operator $\Psi^{*(i)}$
$$%\beq\label{normaliz-dual}
\Phi^{*(i)}(\z)(u_i\ot v_{\ep_{1-i}})=(-\z)^{-i/2}u_{{1-i}}+\cdots,\quad
\ep_0=-,\ \ \ep_1=+\  .
$$%\eeq

It is clear that solution to the equations
\r{fI1}--\r{fI3} can be found in terms of exponential
functions of bosons.
Equation \r{fI1} can be satisfied by the exponent
$$
\exp\left(\sum_{n=1}^{\infty}{a_{-n}\over n}(\z+\h)^n\right)
\ee^{a_0/2} $$
while \r{fI2}--\r{fI3} have a solution   \r{eta} which can be rewritten
using Stirling formula
in the following form
\beq
\stackreb{\mbox{\rm lim}}{N\to \infty}
(2\h)^{\da/2}
\left( {\Ga{\fract{\h-\z}{2\h}}\over\Ga{\fract{-\z}{2\h}}}    \right)^{\da}
\prod_{k=0}^{N}
\exp\left(-\sum_{n=1}^{\infty}  {a_n\over n}
\left[(\z-2k\h)^{-n}-(\z-\h-2k\h)^{-n}\right]    \right)\ .
\label{***}
\eeq
For example,
in order to verify that operator \r{Phi-} is a solution to the equations
\r{fI3} one should substitute \r{***} into l.h.s. of
\r{fI3} and take the limit $N\to\infty$ using the following property
of the ratio of $\Gamma$-functions:
\begin{equation}\label{nice-formula}
\stackreb{\mbox{\rm lim}}{N\to \infty}
(N+b)^{-1}
{ \Gamma^2\left(N+a+\fract{1}{2}\right)
\over
\Gamma^2\left(N+a\right) }   = 1, \quad \forall a,b\in\CC \ \ \mbox{finite}
\ .\end{equation}

Note also that to verify of
the normalization conditions \r{normaliz}
one should expand the ratio of $\Gamma$-functions
in \r{Phi-} in series of positive powers
of $\hbar$ according to the Stirling formula
\begin{equation}\label{Stirling}
 \Gamma(x)\stackreb{=}{x\to\infty}
x^{x-{1\over2}} e^{-x}\sqrt{2\pi}(1+o(1/x))\ .
\end{equation}
\medskip

        \noindent
{\sl Remark.} One can note that equations
\r{fI2}--\r{fI3} have two different solutions. The second one
is
\bea
&\exp\left(\sum_{n=1}^{\infty}{a_{-n}\over n}(\z+\h)^n\right)
\ee^{a_0/2}
\stackreb{\mbox{\rm lim}}{N\to \infty}
(-2\h)^{\da/2}
\left( {\Ga{\fract{\z+2\h}{2\h}}\over\Ga{\fract{\z+\h}{2\h}}} \right)^{\da}
\nn\\
&\quad\times\prod_{k=0}^{N}
\exp\left(-\sum_{n=1}^{\infty}  {a_n\over n}
\left[(\z+\h+2k\h)^{-n}-(\z+2\h+2k\h)^{-n}\right]    \right),
\label{4.19}
\eea
which produces different from \r{fusionI} rule of normal ordering, but the
same commutation relation \r{ZFI}.
To fix a unique solution we should add to the defining equations
\r{fI1}--\r{fI3} some more information. This will be a correspondence
with analogous bosonization formulas for the intertwining operators of
the highest weight modules for the
quantum affine algebra \cite{JM}.
A limit to rational (Yangian)
case from trigonometric (quantum affine algebra)
one means the replacement
$$%\begin{equation}\label{replac}
\w=
 q^{ 4\left({\z\over 2\h}\right)}
$$%\end{equation}
and then sending $q^4\to 1$.    In this limit
$$%\begin{equation}\label{lim}
  \prod_{n=0}^{\infty} (1-q^{4(x+n)}) \sim 
(1-q^4)^{1-x}
\Gamma^{-1}(x)
  \prod_{n=0}^{\infty} (1-q^{4(1+n)})\ .
$$%\end{equation}
The normal ordering relation of
type I intertwining
 operators for quantum affine algebra $\Uq$ produces a factor
$$%\begin{equation}\label{factor-trig}
   \prod_{n=0}^{\infty} {(1- q^{2+4n}\w_1/\w_2)\over (1- q^{4+4n}\w_1/\w_2)}
$$%\end{equation}
which  goes to a necessary ratio of $\Gamma$-functions
\r{fusionI} up to some divergent factor in the limit $q^4\to 1$. 
The choice of the solution \r{4.19} corresponds to the choice of the
intertwining operators for the quantum affine algebra which have
a ``good" analytical properties in the limit $q\to\infty$.
\medskip

Note also that in the classical limit $\h\to0$, our intertwining operators
coincide with those of the affine algebra $\widehat{\frak{sl}}_2$ in
homogeneous Frenkel-Kac realization because of the sort of
the Stirling formula
$$%\beq
\stackreb{\mbox{\rm lim}}{\hbar\to0} \sum_{k=0}^\infty \left(
{1\over \z-2k\h}-{1\over \z-\h-2k\h}\right)={1\over 2\z}\ .
%\label{Stirling2}
$$%\eeq

\setcounter{equation}{0}

\section{Commutation Relations for
                  the Intertwining Ope\-ra\-tors}

Let us formulate the following
\medskip

\noindent
{\bf Proposition 2.} {\it
Type I and type II intertwining operators satisfy following commutation
and normalization relations      }
\begin{eqnarray}
\Phi_{\ep_2}(\z_2) \Phi_{\ep_1}(\z_1)
&=&R_{\ep_1\ep_2}^{\ep'_1\ep'_2} (\z_1-\z_2)
\Phi_{\ep'_1}(\z_1) \Phi_{\ep'_2}(\z_2)\ , \label{ZFI} \\
\Psi^*_{\ep_1}(\b_1) \Psi^*_{\ep_2}(\b_2)
&=&-R_{\ep_1\ep_2}^{\ep'_1\ep'_2} (\b_1-\b_2)
\Psi^*_{\ep'_2}(\b_2) \Psi^*_{\ep'_1}(\b_1)\ , \label{ZFII} \\
\Phi_{\ep_1}(\z_1) \Psi^*_{\ep_2}(\b_2)
&=&\tau(\z_1-\b_2)
\Psi^*_{\ep_2}(\b_2) \Phi_{\ep_1}(\z_1)\ , \label{ZFI-II}\\
g\sum_{\ep}\Phi^{*(1-i)}_\ep(\z)\Phi^{(i)}_\ep(\z)&=&\id
\ ,\label{inverse}\\
\label{orthI}
g\Phi^{(1-i)}_{\ep_1}(\z)\Phi^{*(i)}_{\ep_2}(\z)
&=&\delta_{\ep_1\ep_2}\ \id\ ,\\
\label{orthII}
g^{-1}\Psi^{(1-i)}_{\ep_1}(\b_1)\Psi^{*(i)}_{\ep_2}(\b_2)
&=&{\delta_{\ep_1\ep_2}
\over \b_1-\b_2} + o(\b_1-\b_2)\ ,
\end{eqnarray}
{\it where the $R$-matrix is given by }
\beq\label{ver-R-mat}
R(z)=r(z)\overline R(z)
\eeq
{\it and}
\begin{eqnarray}
r(z)={\Ga{\fract{1}{2}-\fract{z}{2\h}}\Ga{1+\fract{z}{2\h}}\over
        \Ga{\fract{1}{2}+\fract{z}{2\h}}\Ga{1-\fract{z}{2\h}}
   }\ , \qquad
\overline R(z)= \left(
\begin{array}{cccc}
1&0&0&0\\  0&b(z)&c(z)&0\\  0&c(z)&b(z)&0\\  0&0&0&1
\end{array}           \right)      \ ,     \label{R-mat}
\end{eqnarray}
$$ \tau(z)= -\mbox{\rm ctg}\ {\pi z\over2\h}, \quad
b(z)={z\over z+\h},\quad c(z)\ =\  {\h\over z+\h}\ ,
\quad g\ =\ \sqrt{2\h\over\pi}\ .
$$
{\it
The commutation relations between $\Phi^{*(i)}(\z)$ and $\Psi^{(i)}(\b)$
follows from \mbox{\rm \r{ZFI}--\r{ZFI-II}} and the identifications
 \mbox{\rm \r{comparI}} and \mbox{\rm \r{comparII}}.}
\medskip

One can easily check that $R$-matrix \r{ver-R-mat}
satisfy the unitary and crossing symmetry conditions
$$%\beq\label{unitar}
R(z)R(-z)=1,
\quad
(C\ot\id)\, R(z)\, (C\ot\id) =R^{t_1}(-z-h)
$$%\eeq
with charge conjugation matrix
\beq\label{charge}
C=\left(\begin{array}{cc} 0&-1\\1&0 \end{array}\right)
\eeq
and $R^{t_1}(z)$ means  transposition with respect to the first space.

The proof of Proposition 2 is based on  the
normal ordering relations
\bea
\Phi_-(\z_2)\Phi_-(\z_1)      &=&
(2\h)^{1/2}
{\Ga{1+\fract{\z_1-\z_2}{2\h}}\over
\Ga{\fract{1}{2}+\fract{\z_1-\z_2}{2\h}}}
{:}\Phi_-(\z_2)\Phi_-(\z_1){:}\ ,
\label{fusionI}                        \\
\Psi^*_-(\b_1)\Psi^*_-(\b_2)      &=&
(2\h)^{1/2}{\Ga{\fract{1}{2}-\fract{\b_1-\b_2}{2\h}}\over
\Ga{-\fract{\b_1-\b_2}{2\h}}}{:}\Psi^*_-(\b_1)\Psi^*_-(\b_2){:}\ ,
\label{fusionII}                    \\
\Psi^*_-(\b_2)\Phi_-(\z_1)      &=&
(2\h)^{-1/2}{\Ga{\fract{1}{2}+\fract{\z_1-\b_2}{2\h}}\over
\Ga{1+\fract{\z_1-\b_2}{2\h}}}{:}\Psi^*_-(\b_2)\Phi_-(\z_1){:}\ ,
\label{fusionI-II}                    \\
\Phi_-(\z)f(\u)   &=&{1\over \u-\z}\    {:}\Phi_-(\z)f(\u){:}
\ = \
-\sum_{k=0}^{\infty} {\u^k\over \z^{k+1}}
{:}\Phi_-(\z)f(\u){:}\ ,
\nn\\
f(\u)\Phi_-(\z)   &=&{1\over \u-\z-\h}\ {:}\Phi_-(\z)f(\u){:}
\ =  \
\sum_{k=0}^{\infty} {(\z+\h)^k\over \u^{k+1}}
{:}\Phi_-(\z)f(\u){:}   \ ,
\label{norderI}\\
\Psi^*_-(\b)e(\v) &=&{1\over \v-\b-\h}\ {:}\Psi^*_-(\b)e(\v){:}
\ = \
-\sum_{k=0}^{\infty} {\v^k\over (\b+\h)^{k+1}}
{:}\Psi^*_-(\b)e(\v){:}\ ,\nn\\
e(\v)\Psi^*_-(\b)  &=&{1\over \v-\b}\    {:}\Psi^*_-(\b)e(\v){:}
\ =  \
\sum_{k=0}^{\infty} {\b^k\over \v^{k+1}}
{:}\Psi^*_-(\b)e(\v){:} \ ,
\label{norderII}\\
e(\v)\Phi_-(\z)  &=& \Phi_-(\z) e(\v)
 = (\v-\z-\h)\ {:}\Phi_-(\z) e(\v){:}\ ,\nn\\
f(\u)\Psi^*_-(\b)&=& \Psi^*_-(\b)f(\u)=
(\u-\b)\    {:} \Psi^*_-(\b)f(\u){:}\ .\nn
\end{eqnarray}
For example, to check three terms relations
like
$$\Phi_+(\z_2) \Phi_-(\z_1)= r(\z)\left[
c(\z)\Phi_+(\z_1) \Phi_-(\z_2)+b(\z)\Phi_-(\z_1) \Phi_+(\z_2)\right]$$
and
$$\Phi_+(\z_2) \Phi_+(\z_1)= r(\z) \Phi_+(\z_1) \Phi_+(\z_2) \qquad
(\z=\z_1-\z_2)\ $$
we have to  insert  there
an operator
$$\Phi_+(\z,\u)=\Phi_-(\z)f(\u)-f(\u)\Phi_-(\z)$$
instead of the operator $\Phi_+(\z)$. Then using normal ordering
relations  \r{norderI}
we can easily see that the relations in question will be
satisfied as formal series identities with respect to all spectral
parameters $\u$, $\z_1$ and $\z_2$.
Commutativity of type I and type II intertwining operators to the scalar
factor $\tau(\z_1-\z_2)$  \r{ZFI-II} follows from the commutativities
of the operators $\Phi^{(i)}_-(\z)$ and $\Psi^{*(i)}_-(\b)$ with
currents $e(\u)$ and $f(\v)$ respectively.
The proof of
the identities \r{inverse} and \r{orthI} is given
in the Appendix.

\setcounter{equation}{0}
\section{Universal ${\cal R}$-Matrix Formulation of $\DYsect$}

The goal of this section is to describe
$\Yd$ in ``$RLL$'' formalism      \cite{FRT,RS}.
The essential part of this construction is the notion of the
universal ${\cal R}$-matrix.
The universal ${\cal R}$-matrix for the central extended Yangian
double have been obtained in \cite{KT,K} by inves\-ti\-gat\-ing
the canonical pairing in $\Yd$.
This ${\cal R}$-matrix is
\beq
{\cal R} = {\cal R}_+\cdot {\cal R}_0 \cdot \exp(\h c\ot d) \cdot {\cal R}_-
\ ,\label{univ-R-matr}
\eeq
where
\bea
{\cal R}_+&=& \prod_{k\geq 0}^\rightarrow \exp (-\h e_k\ot f_{-k-1})
= \exp(-\h e_0\ot f_{-1})\exp(-\h e_1\ot f_{-2})\cdots
\ ,\nn\\%\label{R_+}
{\cal R}_-&=& \prod_{k\geq 0}^\leftarrow \exp (-\h f_k\ot e_{-k-1})
= \cdots \exp(-\h f_1\ot e_{-2})\exp(-\h f_0\ot e_{-1})
\ ,\nn\\%\label{R_-}
{\cal R}_0&=& \prod_{n\geq 0} \exp
\left(\Res{u=v}
\left( {d\over du} \ln\, h^+(u)\ot
\ln\,h^-(v+\h+2n\h)\right)\right)
\ .\nn%\label{R_0}
\eea
Here a residue operation $\Res{}$ is defined as
$$\Res{u=v}\left(\sum_{i\geq0} a_i u^{-i-1}\ot
\sum_{k\geq0} b_k v^{k}
\right) =
\sum_{i\geq0} a_i\ot b_i\ . $$

In the evaluation representation $\pi(z)$ of $\Yd$ in the space $V_z$
the operators $e_k$, $f_k$, $h_k$, $c$, $k\in\ZZ$
act as follows
$$
e_k v_{+,z}=f_k v_{-,z}=c v_{\pm,z}=0,\quad
e_k v_{-,z}=z^k v_{+,z},\quad f_k v_{+,z}=z^k v_{-,z},\quad
h_k v_{\pm,z}=\pm z^k v_{\pm,z}\ .
$$
With this action the $L$-operators
(see \cite{FR} for $U_q(\widehat{\frak{g}})$ case)
\beq\label{7.00}
L^-(z)= (\pi(z)\ot \id) {\cal R} \exp(-\h c\ot d),\qquad
L^+(z)= (\pi(z)\ot \id) \exp(\h d\ot c) \left({\cal R}^{21}\right)^{-1}
\eeq
appear in the Gauss decomposed form
\bea
L^-(z)&=& \left(\begin{array}{cc} 1& \h f^-(z)\\ 0&1\end{array}\right)
\left(\begin{array}{cc} (k^-(z-\h))^{-1}&0\\ 0& k^-(z)
\end{array}\right)
\left(\begin{array}{cc} 1&0\\ \h e^-(z)&1\end{array}\right)\ ,
\label{GaussL-}\\
L^+(z)&=&
\left(\begin{array}{cc} 1&\h f^+(z-c\h)\\ 0&1\end{array}\right)
\left(\begin{array}{cc} (k^+(z-\h))^{-1}&0\\ 0& k^+(z)
\end{array}\right)
\left(\begin{array}{cc} 1&0\\ \h e^+(z)&1\end{array}\right)\ .
\label{GaussL+}
\eea
To obtain the diagonal part of the operators $L^+(z)$ the
property of operation Res have been used
$$
\Res{u=v}\left({d\over du}
\sum_{i\geq0} a_i u^{-i-1}\ot
\sum_{k\geq0} b_k v^{k}
\right) = -
\Res{u=v}\left(\sum_{i\geq0} a_i u^{-i-1}\ot
{d\over dv}\sum_{k\geq0} b_k v^{k}
\right) \ .
$$
From formulas \r{GaussL-} and \r{GaussL+}
the following representation  of
the fields $k^\pm(z)$ can be obtained
\beq\label{k-fields}
k^-(z)=\prod_{n=0}^\infty {h^-(z+2\h+2n\h)\over h^-(z+\h+2n\h)},\qquad
k^+(z)=\prod_{n=0}^\infty {h^+(z-\h-2n\h)\over h^+(z-2n\h)}\ .
\eeq
To proceed further we have to calculate the universal ${\cal R}$-matrix
on the tensor product of two evaluation representations ($z=z_1-z_2$)
\bea
\tilde R^-(z)&=&(\pi(z_1)\ot \pi(z_2)){\cal R}\exp(-\h c\ot d)\ = \
(\id\ot \rho(z_2)) L^-(z_1)\ =\ \rho^-(z)\overline R(z)
\nn%\label{R-}
\\
\tilde R^+(z)&=&(\pi(z_1)\ot
\pi(z_2))\exp(\h d\ot c)\left({\cal R}^{21}\right)^{-1}
\ =\ (\id\ot \rho(z_2)) L^+(z_1)\ =\ \rho^+(z)\overline R(z)
\ ,\nn%\label{R+}
\eea
where   $\overline R(z)$ is given by   \r{R-mat} and
$$%\begin{equation}\label{univR-norm}
      \rho^\pm(z)=
\left[{\Ga{\mp\fract{z}{2\h}}\Ga{1\mp\fract{z}{2\h}}\over
\Gamma^2\left(\fract{1}{2}\mp\fract{z}{2\h}\right)}\right]^{\mp1}\ .
$$%\end{equation}

For the $R$-matrices $\tilde R^\pm(z)$ the conditions of unitary
and crossing symmetry look as follows
\beq
R^+(z) R^-(-z)=1,
\qquad(C\ot\id)\, \tilde R^\pm\, (C\ot\id)^{-1}=
\tilde R^{\mp\, t_1}(-z-\h)
\label{cross-univ}
\eeq
with charge conjugation matrix $C$ defined by \r{charge}.

The following properties of the universal $R$-matrix ${\cal R}$ and
representations of $\Yd$
\begin{itemize}
\item[(i)] Yang-Baxter relation for ${\cal R}$,
\item[(ii)]  comultiplication     rules for ${\cal R}$
$$(\D\otimes \id)\R =\R^{13}\R^{23}, \hsp
        (\id\otimes\D)\R =\R^{13}\R^{12}\ ,$$
\item[(iii)] realization of the trivial representation of $\Yd$ as
submodule
$${\CC}\ (v_+\ot v_--v_-\ot v_+)\hookrightarrow \pi(z+\h)\ot\pi(z)$$
\end{itemize}
imply the following relations for $L^\pm(z)$:
\begin{eqnarray}
\tilde R^-(z+c\h)L^-_1(z_1)L^+_2(z_2)&=& L^+_2(z_2)L^-_1(z_1) \tilde R^-(z)\ ,
\label{RL-L+}
\\
\tilde R^\pm(z)L^\pm_1(z_1)L^\pm_2(z_2)&=& L^\pm_2(z_2)L^\pm_1(z_1)
\tilde R^\pm(z) \  ,
\label{RLL}    \\
\mbox{q-det}\,L^\pm(z)&=&1\ .
\label{q-detL}
\end{eqnarray}
Here $\mbox{q-det}\,A(z)$ is equal to
$A_{11}(z+\h)A_{22}(z)-A_{12}(z+\h)A_{21}(z)$ and
comultiplication rules are
\bea
\Delta' L^+(z)&=&L^+(z-\h c_2)
\dot\otimes L^+(z)\ , \quad
 \Delta l_{ij}^+(z)\ =\
\sum_kl^+_{kj}(z)\otimes l_{ik}^+(z-\h c_1)\ ,
\label{D-L+}     \\
\Delta' L^-(z)&=&L^-(z)\dot\otimes L^-(z)\ , \quad
 \Delta l_{ij}^-(z)\ =\ \sum_kl_{kj}^-(z)\otimes l_{ik}^-(z)\ ,
\label{D-L-}
\eea
where $\Delta'$ means the flipped comultiplication and $\dot\ot$ is
matrix tensor product. Using Gauss decomposition of $L$-operators
we can easily obtain from \r{D-L+} and \r{D-L-} the full
comultiplication formulas for the currents $e^\pm(u)$, $f^\pm(u)$,
$h^\pm(u)$ (see \r{3.6}).

Conversly, we can start from an algebra of coefficients of the operators
$L^\pm(z)$ subjected to \r{RL-L+}--\r{q-detL} and prove that the currents
\begin{eqnarray}
e(z)&=&e^+(z)-e^-(z)\ ,\nn\\
f(z)&=&f^+(z)-f^-(z)\ ,\nn\\
h^\pm(z)&=&[k^\pm(z)k^\pm(z-h)]^{-1}
\label{Frenk-Ding}
\end{eqnarray}
defined by the Gauss decompositions
\r{GaussL-}, \r{GaussL+} satisfy the relations
of $\Yd$. These calculations (Yangian analogs of Ding--Frenkel
arguments \cite{DF}) are presented in the Appendix.

\setcounter{equation}{0}
\section{Miki's Formulas and Bosonization of $L$-operators}

Following the approach of the paper \cite{M} let us compose from the
intertwining operators the following $2\times2$ operator
End$V_i$-valued matrices
\beq
L^+_{\ep\nu}(z)=g\Phi_\ep(z-\h)\Psi^*_\nu(z), \qquad
L^-_{\ep\nu}(z)=g\Psi^*_\nu(z)\Phi_\ep(z)\ . \label{Miki}
\eeq
One can check from commutation and orthogonality relations
\r{ZFI}--\r{ZFI-II} that these operators  satisfy the relations
\r{RL-L+}--\r{q-detL} with $c=1$.

It is important to point out here that the order of intertwining operators
in \r{Miki} is strictly fixed, because only in this order
the normal ordering of these operators evaluated in the same points
is well defined. It follows
from the fusion rules
\begin{eqnarray}
\Phi_-(z_1-\h)\Psi^*_-(z_2)&=& (2\h)^{-1/2}
{\Ga{\fract{1}{2}-\fract{z_1-z_2}{2\h}}\over\Ga{1-\fract{z_1-z_2}{2\h}}}
{:}\Phi_-(z_1-\h)\Psi^*_-(z_2){:} \ ,    \nn\\
\Psi^*_-(z_1)\Phi_-(z_2)&=& (2\h)^{-1/2}
{\Ga{\fract{1}{2}-\fract{z_1-z_2}{2\h}}\over\Ga{1-\fract{z_1-z_2}{2\h}}}
{:}\Psi^*_-(z_1)\Phi_-(z_2-){:}  \ .   \nn%  \label{fusI-II}
\end{eqnarray}
Also it is necessary to say that checking \r{RLL} or \r{RL-L+}
by means of \r{Miki} one will never interchange the intertwining operators
evaluated in the same points which is an ill defined operation.

Our goal now is to show that $L$-operators \r{Miki} coincide with those
defined by the universal ${\cal R}$-matrix \r{univ-R-matr} and evaluated
at level one modules $V_0$ and $V_1$. In order to do this we express
operators \r{Miki} in terms of free bosons using \r{Phi-} and \r{Psi*-}.
First calculate operators $L_{--}^\pm(z)$.
\begin{eqnarray}
L^+_{--}(z)&=&g\Phi_-(z-\h)\Psi^*_-(z)={\eta_+(z)\over\eta_+(z-\h)}\ ,\nn\\
L^-_{--}(z)&=& g\Psi^*_-(z)\Phi_-(z) =  {\phi_-(z+\h)\over \phi_-(z)}\ .
\label{k-}
\end{eqnarray}
Comparing \r{k-} with \r{k-fields} and \r{h-fields} we see that
$$
L^+_{--}(z)=k^+(z),\qquad
L^-_{--}(z)=k^-(z),\qquad (c=1)\ .
$$

To write down the bosonized expressions for the rest elements of the
$L$-operators we first have to obtain more reasonable for our purposes
representation of the ``$+$'' components of the intertwining
operators $\Phi(z)$ and $\Psi^*(z)$.
Using definition \r{inteq}
for $x=f^\pm(z)$, $e^\pm(z)$, explicit formulas of comultiplication
\r{3.6} and \r{eval} we can obtain following
relations
\bea
{1\over z_2-z_1}\Phi_+(z_1)&=&\Phi_-(z_1)f^\pm(z_2)-
{z_2-z_1-\h\over z_2-z_1} f^\pm(z_2)\Phi_-(z_1)\ ,
\quad \vert z_2\vert\probaF \vert z_1  \vert
\label{6.4a}
\\
{1\over z_2-z_1}\Psi^*_+(z_1)&=&e^\pm(z_2)\Psi^*_-(z_1)-
{z_2-z_1-\h\over z_2-z_1} \Psi^*_-(z)e^\pm(z_2)\ ,
\quad \vert z_2\vert\probaF \vert z_1  \vert
\label{6.4b}
\\
0&=&\Phi_-(z_1)e^\pm(z_2)- e^\pm(z_2)\Phi_-(z_1)\ , \label{6.4c}
\\
0&=&f^\pm(z_2)\Psi^*_-(z_1)- \Psi^*_-(z)f^\pm(z_2)\ ,  \label{6.4d}
\eea
where $f^\pm(z_2)$ and $e^\pm(z_2)$ defined by
\r{2.0}.

Let us comment  the meaning of the formulas
\r{6.4a}, for example.
The formulas \r{6.4b} can be treated analogously. Expanding
first equation in
\r{6.4a} in positive powers of $z_1/z_2$ and second one in
\r{6.4a} in positive powers of $z_2/z_1$ we obtain in addition
to \r{4.16}  the relations
\bea
\Phi_+(z_1)&=&[\Phi_-(z_1),f_{m}]z_1^{-m}+
\h\sum_{k=0}^{m-1}z_1^{-k-1}f_{k}\Phi_-(z_1)\ ,
\nn \\
\Phi_+(z_1)&=&[\Phi_-(z_1),f_{-m}]z_1^{m}-
\h\sum_{k=0}^{m-1}z_1^{k}f_{-k-1}\Phi_-(z_1)\ ,
\quad m=1,2,\ldots
\nn
\eea
These formulas define infinitly many relations between the operators
$\Phi_-(z_1)$ and $f_{m}$, $m\in\ZZ$ which can be encoded in a single
relation
between      the operator
$\Phi_-(z_1)$ and the total current $f(u)$
\beq
\label{*}
(z_2-z_1) \Phi_-(z_1) f(z_2)=(z_2-z_1-\h)f(z_2)\Phi_-(z_1) .
\eeq
This relation is in an agreement with normal ordering rules
\r{norderI} and therefore the components of the intertwining
operators
\r{Phi-}--\r{comparII} satisfy the intertwining relations for
all the generators of the Yangian double
$\Yd$.

Using now the relations
\r{6.4a}--\r{6.4d} and also
\r{fI1}, \r{fI2}
  after
 straightforward calculations one can obtain
the rest elements of $L$-operators
\bea
L^+_{-+}(z)&=&g\Phi_-(z-\h)\Psi^*_+(z)=\h k^+(z) e^+(z)\ ,\nn\\
L^+_{+-}(z)&=&g\Phi_+(z-\h)\Psi^*_-(z)=\h f^+(z-\h)k^+(z)\ ,\nn\\
L^-_{-+}(z)&=&g\Psi^*_+(z)\Phi_-(z)=\h k^-(z) e^-(z)\ ,\nn\\
L^-_{+-}(z)&=&g\Psi^*_-(z)\Phi_+(z)=\h f^-(z)k^-(z)\ ,\nn\\
L^+_{++}(z)&=&g\Phi_+(z-\h)\Psi^*_+(z)
=\left(k^+(z-\h)\right)^{-1}+ \h^2 f^+(z-\h)k^+(z)e^+(z)\ ,\nn\\
L^-_{++}(z)&=&g\Psi^*_+(z)\Phi_+(z)
=\left(k^-(z-\h)\right)^{-1}+ \h^2 f^-(z)k^-(z)e^-(z)\ ,\nn
\eea
where operator identities \r{Frenk-Ding} have been used.
From these formulas follow that the Gauss decomposition
of the $L^\pm(z)$-operators \r{Miki} coincides with those in
\r{GaussL-} and \r{GaussL+} at the level $c=1$.

Let us comment the using of the relations \r{6.4a} and \r{6.4b}.
First of all multiply them by $z_2-z_1$
to obtain\footnote{In the previous version of this paper we used
incorrect formulas
\bea
\Phi_+(z)&=&\h f^+(z)\Phi_-(z)=\h \Phi_-(z)f^-(z+\h)\ ,
\nn  \\
\Psi^*_+(z)&=&\h \Psi^*_-(z) e^+(z)=\h e^-(z+\h)\Psi^*_-(z)\ ,
\nn
\eea
from which one cannot obtain the correct formulas for
elements of $L$-operators $L^\pm_{++}(z)$. }
\bea
\Phi_+(z_1)
&=&(z_2-z_1)\Phi_-(z_1)f^\pm(z_2)-
(z_2-z_1-\h) f^\pm(z_2)\Phi_-(z_1)\ ,
\label{6.4am}
\\
&\equiv&\Phi_-(z_1)f_0 -f_0\Phi_-(z_1)\ ,
\label{6.4an}
\\
\Psi^*_+(z_1)
&=&(z_2-z_1)e^\pm(z_2)\Psi^*_-(z_1)-
(z_2-z_1-\h) \Psi^*_-(z_1)e^\pm(z_2)\ , \label{6.4bm}
\\
&\equiv&e_0\Psi^*_-(z_1)- \Psi^*_-(z_1)e_0\ .
\label{6.4bn}
\eea
The whole set of the identities
\r{*} ensures the equations \r{6.4an}
and \r{6.4bn}.
Let us substitute operator $\Phi_+(z_1)$
in the form
\r{6.4am}
into, for example,
 $L^-_{+-}(z_1)$. We obtain
\bea
L^-_{+-}(z_1)&=&g\Psi^*_-(z_1)\Phi_+(z_1) \nn\\
           &=&(z_2-z_1)k^-(z_1)f^-(z_2)-
(z_2-z_1-\h) f^-(z_2)k^-(z_1)=\h f^-(z_1)k^-(z_1)\ ,\nn
\eea
where identity \r{com5} have been used.

We can justify above formulas for comultiplication by the Miki's relations
for the $L$-operators at level 1. Indeed, using also \r{Miki},
\r{ZFI}--\r{ZFI-II} and identities
$$ r(z)\tau(-z)=\rho^-(z)
\qquad\mbox{and} \qquad
r(z)\tau(-z-\h)=\rho^+(z)$$
we can calculate
\begin{eqnarray}
\Phi_{\ep_2}(z_2)L^-_{\ep_1\nu}  (z_1)
&=&
\tilde R_{\phantom{-\,}\ep_1\ep_2}^{-\,\ep'_1\ep'_2}(z_1-z_2)
L^-_{\ep'_1\nu}(z_1)\Phi_{\ep'_2}(z_2)\ ,\label{L-Phi} \\
\Phi_{\ep_2}(z_2)L^+_{\ep_1\nu}  (z_1)
&=&
\tilde R_{\phantom{+\,}\ep_1\ep_2}^{+\,\ep'_1\ep'_2}(z_1-z_2-\h)
L^+_{\ep'_1\nu}(z_1)\Phi_{\ep'_2}(z_2)\label{L+Phi}
\end{eqnarray}
which obviously coincide with those formulas which can be obtained
by the definition of the intertwining operators \r{inteq}
and \r{D-L+} and \r{D-L-}.

\setcounter{equation}{0}

\section{Conclusion}
To conclude let us shortly repeat the main result obtained in this paper.
Starting from the free field
representation of the central extension of the Yangian double at level
$c=1$ \cite{K} the explicit formulas for the operators which intertwine these
representations have been obtained. As well as in the case of quantum
affine algebras these operators appear to be divided in two types.
It was shown that type II intertwining operators satisfy the commutation
relations of the Zamolodchikov-Faddeev algebra for $S$-matrix
coincided  with kink scattering in $SU(2)$-invariant Thirring
model.
$$%\begin{equation}\label{S-mat}
S_{12}(\beta)=
{\Ga{\fract{1}{2}+\fract{\beta}{2\pi i}}\Ga{-\fract{\beta}{2\pi i}}\over
        \Ga{\fract{1}{2}-\fract{\beta}{2\pi i}}\Ga{\fract{\beta}{2\pi i}}    }
\left[{\beta-\pi i P_{12}\over \beta-\pi i}\right]  \ ,
$$
where the deformation parameter $\h$ of $\Yd$ should be chosen
as
$$\h=-\pi i\ .$$
We can state (as well as it was done in \cite{JM} for lattice
XXZ model) that the physical interesting quantities of $SU(2)$-invariant
Thirring  model,
such as form-factors of the local operators, can be calculated
as traces over infinite dimensional representation spaces of $\Yd$
of the weighted product of the types I and II intertwining operators.
We address the calculation of these traces to the forthcoming
paper \cite{KLP}.

\section{Acknowledgments}

Part of this work was done while first two authors were visiting
Roma I University in September--October 1995.
These authors would like to acknowledge Prof. A. Degasperis
for the hospitality during this visit.
The first author (S.K.) was supported by RFBR grant 95-01-01106
and INTAS grant 93-0023.
The second author (D.L.)
 was suppoted by  RFBR grant 95-01-01101, INTAS
grant 93-0166, Wolkswagen Stiffung and AMSFSU grant.
The third author (S.P.) would like to acknowledge Departamento de F\'\i sica
Te\'orica of Universidad de Zaragoza where part of this work was done.
The research of S.P. was supported by Direcci\'on General de
Investigaci\'on Cient\'\i fica y T\'ecnica (Madrid)
and INTAS grant 93-2058.
\bigskip
\medskip

\noindent
{\it Note added:
Before this paper was resubmitted to the electronic archive
the authors
found that the central extension of the Yangian double
associated with algebra $\widehat{\frak{gl}}_2$
and free field
representation of this algebra at level one have been obtained
independently by
K.~Iohara and H.~Kohno \cite{IK}
using the methods different  from ones in \cite{K}.
We thank Prof. M.~Jimbo for pointing out our attention to this
publication.}
\bigskip
\medskip

\setcounter{equation}{0}
\renewcommand{\theequation}{\mbox{\rm A.{\arabic{equation}}}}

\app{}
\sapp{Orthogonality relations for intertwining operators}

\noindent
Let us prove \r{orthI} in the simplest case $-\ep_1=\ep_2=+$.
Because of the relation
\bea
\Phi^{(1-i)}_-(z_2)\Phi^{*(i)}_+(z_1)&=&
(-1)^i\Phi_-(z_2)\Phi_-(z_1-\h)\nn\\
&=&
(-1)^{i}
(2\h)^{-1/2}
(z_1-z_2){\Ga{\fract{1}{2}+\fract{z_1-z_2}{2\h}}\over
\Ga{1+\fract{z_1-z_2}{2\h}}}
{:}\Phi_-(z_2)\Phi_-(z_1-\h){:}
\label{ort1}
\eea
which follows from \r{comparI} and \r{fusionI}
we have
$$\Phi^{(1-i)}_-(z)\Phi^{*(i)}_+(z) = 0$$
if we suppose that normal ordered product ${:}\Phi_-(z_2)\Phi_-(z_1-\h){:}$
cannot produces  pole in $z_1-z_2$ while acting in $V_0$ or $V_1$.
Here we have to be more accurate. In fact we need the following
lemma.

\medskip
\noindent
{\bf Lemma.} {\it If two operators $A(z_1)$ and $B(z_2)$
satisfy the relation
$$A(z_1)B(z_2)=(z_1-z_2)\ {:}A(z_1)B(z_2){:}$$
then the product $A(z)B(z)$ acts on any finite element of the representation
spaces $V_0$ and $V_1$ multiplying it by zero
$$A(z)B(z) v =0,\quad \forall v \in V_0\quad\mbox{or}\quad V_1\ .$$}

\noindent
The word {\it finite} in the formulation of the Lemma means
finite sum of the elements of the type \r{modules}.

To prove \r{orthI} in more complicated case $\ep_1=\ep_2=+$ we
have to substitute operator
$$\Phi^{(1-i)}
_+(z_2)=\res{v}\Phi^{(1-i)}_+(z_2,u)=\res{v}
\left[\Phi^{(1-i)}_-(z_2)f(u)-f(u)\Phi^{(1-i)}_-(z_2)\right]$$
into product $\Phi^{(1-i)}_+(z_2)\Phi^{*(i)}_+(z_1)$.
Residue operation $\res{v}A(v)$ means
\begin{equation}\label{residue}
\res{v} A(v)=A_0,\qquad\mbox{if}\qquad  A(v) = \sum_{k}A_kv^{-k-1}\ .
\end{equation}
 After
normal ordering all products using  formulas  \r{norderI}, \r{ort1}
we obtain
\bea
g\Phi^{(1-i)}_+(z_2,u)\Phi^{*(i)}_+(z_1)  &=& (-1)^i
{\Ga{\fract{1}{2}+\fract{z_1-z_2}{2\h}}\over
\sqrt{\pi}\Ga{1+\fract{z_1-z_2}{2\h}}}
\left(\delta(u-z_2)+ o(z_1-z_2)  \right)\nn\\
&\times& {:}\Phi_-(z_2)f(u)\Phi_-(z_1-\h){:}\ .   \nn%\label{prove1}
\eea
Because of the $\delta$-function we can set $u=z_2$ in the normal
ordered product of the
 operators, set $z_1=z_2$, use the operator
identity
$$%\begin{equation}\label{identI}
{:}\Phi_-(z)f(z)\Phi_-(z-\h){:}= (-1)^{\da}  \nn
$$%\end{equation}
and apply Lemma 1 to obtain that
\begin{equation}
g\Phi^{(1-i)}_+(z)\Phi^{*(i)}_+(z)  =  1
\end{equation}
because operator $(-1)^{\da}$ act on the spaces $V_i$  multiplying  them
by $(-1)^i$.
Other cases can be proved analogously.

Let us prove now the identity \r{inverse}. Because of the fusion rule
$$
\Phi_-(z-\h)\Phi_-(z)=g^{-1}\h\, {:}\Phi_-(z-\h)\Phi_-(z){:}
$$
we have
\bea
&g\sum_{\ep}\Phi^{*(1-i)}_\ep(z)\Phi^{(i)}_\ep(z)=
g(-1)^i\left( \Phi_-(z-\h)\Phi_+(z)-\Phi_+(z-\h)\Phi_-(z)
\right)
\nn\\
&\quad =\ (-1)^i\h\res{u}\left[
{2\over (u-z-\h)(u-z+\h)}
-{1\over (u-z-\h)(u-z)}-{1\over (u-z+\h)(u-z)}
\right] \nn\\
&\qquad\times\
{:}\Phi_-(z-\h)f(u)\Phi_-(z){:}\nn\\
&\quad =(-1)^i\res{u}\delta(u-z)
{:}\Phi_-(z)f(u)\Phi_-(z-\h){:}=  \id \ .\nn
\eea
\bigskip

\sapp{Ding--Frenkel Equivalence}

\noindent
To prove the Ding--Frenkel equivalence of new and $RLL$ realizations
of $\Yd$
we have to insert $L$-operators in the form
$$%\begin{equation}\label{L+}
L^+(z)=
\left(\begin{array}{cc}
\big(k^+(z-\h)\big)^{-1}+\h^2 f^+(z-\h c)k^+(z)e^+(z)&\h f^+(z-\h c)k^+(z)\\
\h k^+(z)e^+(z)& k^+(z)
\end{array}\right)
$$%\end{equation}
and
$$%\begin{equation}\label{L-}
L^-(z)=
\left(\begin{array}{cc}
\big(k^-(z-\h)\big)^{-1}+\h^2 f^-(z)k^-(z)e^-(z)&\h f^-(z)k^-(z)\\
\h k^-(z)e^-(z)& k^-(z)
\end{array}\right)
$$%\end{equation}
into \r{RLL} and \r{RL-L+}.
As usual, $L_1(z)$ and $L_2(z)$ are tensor products
$L(z)\ot 1$ and $1\ot L(z)$ respectively.

First, we will write the commutation relations between $k^\pm(z)$
which follows from $D_1D_2$ commutation relations.
We have
\bea
\label{com1}
k^\pm(z_1)k^{\pm}(z_2)&=&
k^\pm(z_2)k^{\pm}(z_1)\ ,\\
\label{com2}
\rho(z+c\h)k^-(z_1)k^+(z_2)&=&
k^+(z_2)k^-(z_1) \rho(z)\ .
\eea
Note that non-triviality of the commutation relation in \r{com2} is
dictated by presence of
$c$. If $c=0$ as happen for evaluation representations
then all $k$'s commute.
Now note that
commutation relation between $h^+(z_1)$ and $h^-(z_2)$ easily follows from
commutation relations \r{com1}--\r{com2}
if we identify
\begin{equation}\label{newCartan}
\left(k^\pm(z+\h)k^\pm(z)\right)^{-1} = h^\pm(z)
\end{equation}
and by making use the identity
$$\rho^\pm(z+\h)\rho^\pm(z)b(z)=1.$$

To begin with  nontrivial part of commutation relations we will consider
the relation
$$b(z)D(z_1)B(z_2)+ c(z) B(z_1)D(z_2)=a(z)B(z_2)D(z_1)$$
 written
in three cases of \r{RLL} and \r{RL-L+}. It yields ($z=z_1-z_2$)
\bea
b(z)  k^-(z_1)f^-(z_2)k^-(z_2)+ c(z)
f^-(z_1) k^-(z_1) k^-(z_2) &=&
f^-(z_2) k^-(z_2) k^-(z_1)\ ,
\label{com3} \\
b(z)  k^+(z_1)f^+(z_2-c\h)k^+(z_2)+ c(z)
f^+(z_1-c\h) k^+(z_1) k^+(z_2) &=&
f^+(z_2-c\h) k^+(z_2) k^+(z_1)\ ,
\label{com3a}         \\
 \rho(z+c\h)\left[
b(z+c\h)  k^-(z_1)f^+(z_2-c\h)k^+(z_2)\right.&+&\left. c(z+c\h)
f^-(z_1) k^-(z_1) k^+(z_2) \right]     \nn\\
&=&\rho(z)f^+(z_2-c\h) k^+(z_2) k^-(z_1)\ .
\label{com4}
\eea
Multiplying  \r{com3} and \r{com3a} from the right by
$(k^\pm(z_1)k^\pm(z_2))^{-1}$,
\r{com4} by
$(k^-(z_1))^{-1}(k^+(z_2))^{-1}$  and using in the
latter case the commutation relation \r{com2} to cancel the factor
$\rho(z)$  we obtain
\bea
f^-(z_2)&=&
b(z)  k^-(z_1) f^-(z_2)(k^-(z_1))^{-1}+ c(z) f^-(z_1)\ ,
\label{com5}\\
f^+(z_2-c\h)&=&
b(z)  k^+(z_1) f^+(z_2-c\h)(k^\pm(z_1))^{-1}+ c(z) f^+(z_1-c\h)\ ,
\label{com5a}\\
f^+(z_2-c\h)&= &
b(z+c\h)  k^-(z_1) f^+(z_2-c\h)(k^-(z_1))^{-1}+c(z+c\h)
f^-(z_1)                                               \ .
\label{com6}
\eea
Making comparison of
\r{com5} and \r{com6} (in the latter relation we
have to shift $z_2\to z_2+ c\h$) we can see
that the combination
$$%\begin{equation}
f^+(z)  -
f^-(z)
=  f(z) = \sum_{n\in\ZZ}
f_n z^{-n-1}
%\label{cur+}
$$%\end{equation}
has ``good'', namely, vertex like commutation relation with $(k^-(z))^{-1}$
$$%\begin{equation}\label{com7}
(k^-(z_1))^{-1}f(z_2) = b(z) f(z_2)(k^-(z_1))^{-1}.
$$%\end{equation}
To obtain the commutation relation $k^+(z_1)$ with $f(z_2)$ we have
to consider another $DB$ relation, namely,
$$c(z)D(z_1)B(z_2)+ b(z) B(z_1)D(z_2)=a(z)D(z_2)B(z_1)$$
 in the case
of \r{RL-L+}. It is
\begin{eqnarray}
& \rho(z+c\h)\left[
c(z+c\h)  k^-(z_1)f^+(z_2-c\h)k^+(z_2)+b(z+c\h)
f^-(z_1) k^-(z_1) k^+(z_2) \right]     \nn\\
&\quad\quad\quad\quad=\rho(z)k^+(z_2) f^-(z_1) k^-(z_1)\ .
\label{com8}
\end{eqnarray}
Then inserting in l.h.s. of \r{com8}
$(k^+(z_1)) ^{-1}
k^+(z_1)$ between $f^+(z_2-c\h)$ and $k^+(z_2)$, using \r{com2}
to cancel $\rho(z+c\h)$ and multiplying from the right by
$(k^-(z_1))^{-1}(k^+(z_2))^{-1}$    we obtain
\begin{equation}
\label{com9}
c(z+c\h)  k^-(z_1)f^+(z_2-c\h)(k^-(z_1))^{-1}+b(z+c\h) f^-(z_1) =
k^+(z_2) f^-(z_1)(k^+(z_2))^{-1}\ .
\end{equation}
Expressing the combination $k^-(z_1)f^+(z_2-c\h)(k^-(z_1))^{-1}$
from \r{com6} via $f^+(z_2-c\h)$ and $f^-(z_1)$,
replacing $z_1\to z_2$, $z_2\to z_1$ in \r{com9}, and
using the identities
$$ {b(-z)\over b^2(-z)-c^2(-z)}=b(z),\quad
{-c(-z)\over b^2(-z)-c^2(-z)}=c(z)$$
 we obtain
\begin{equation}
\label{com10}
f^-(z_2)=
b(z-c\h)  k^+(z_1) f^-(z_2)(k^+(z_1))^{-1}+ c(z-c\h) f^+(z_1-c\h)\ .
\end{equation}
From \r{com10} and \r{com5a} (in the latter equation we have to shift $z_2
\to z_2+c\h$) we have
$$%\begin{equation}\label{com11}
(k^+(z_1))^{-1}f(u)(z_2) = b(z-c\h) f(z_2)(k^+(z_1))^{-1}.
$$%\end{equation}
Now fourth and fifth lines in \r{DY2} follow from
definition \r{newCartan} and simple identity
\begin{equation}
\label{ide}
b(z-\h)b(z)={z-\h\over z+\h}\ .
\end{equation}

The relations
\bea
e^\pm(z_2)&=&
b( z)  (k^\pm(z_1))^{-1}e^\pm(z_2)k^\pm(z_1)+
c(z) e^\pm(z_1)\ ,\label{com100}
\\
e^\mp(z_2)&=&
b( z)  (k^\pm(z_1))^{-1}e^\mp(z_2)k^\pm(z_1)+
c(z) e^\pm(z_1)             \nn
\eea
which follow from $DC$ relations
 yield for the combination
$$%\begin{equation}
e^+(z)- e^-(z)=  e(z) =  \sum_{n\in\ZZ}
e_n z^{-n-1}
%\label{cur-}
$$%\end{equation}
the commutation relations with $k^\pm(z)$
$$%\begin{equation}\label{com12}
(k^\pm(z_1))^{-1}e(z_2)=b(z)^{-1}e(z_2)(k^\pm(z_1))^{-1}
$$%\end{equation}
and third line in \r{DY2} is again a consequence
of the identity \r{ide}.

\end{document}